\documentclass[11pt]{amsart}

\usepackage[english]{babel}
\usepackage[utf8]{inputenc}
\usepackage{amsmath}
\usepackage{amssymb}
\usepackage{amsfonts}
\usepackage{amsthm}
\usepackage{mathrsfs}
\usepackage[all]{xy}
\usepackage[pdftex]{graphicx}
\usepackage{color}
\usepackage{cite}
\usepackage{url}
\usepackage{indent first}
\usepackage[labelfont=bf,labelsep=period,justification=raggedright]{caption}
\usepackage[english]{babel}
\usepackage[utf8]{inputenc}
\usepackage{hyperref}
\usepackage[colorinlistoftodos]{todonotes}
\makeatletter
\newtheorem*{rep@theorem}{\rep@title}
\newcommand{\newreptheorem}[2]{%
\newenvironment{rep#1}[1]{%
 \def\rep@title{#2 \ref{##1}}%
 \begin{rep@theorem}}%
 {\end{rep@theorem}}}
\makeatother

\DeclareMathOperator{\bias}{bias}

\DeclareMathOperator{\ext}{ex}

\DeclareMathOperator{\Tr}{Tr}

\DeclareMathOperator{\NP}{\textsf{NP}}

\newcommand{\eps}{\epsilon}

\newcommand{\mcB}{\mathcal{B}}

\newcommand{\mcH}{\mathcal{H}}

\newcommand{\mcR}{\mathcal{R}}

\renewcommand{\P}{\mathbf{P}}
\makeatletter
\newcommand\thankssymb[1]{\textsuperscript{\@fnsymbol{#1}}}
\makeatother

\theoremstyle{plain}
\newtheorem{thm}{Theorem}
\newreptheorem{theorem}{Theorem}
\newtheorem{lemma}[thm]{Lemma}
\newtheorem{lem}[thm]{Lemma}
\newtheorem{cor}[thm]{Corollary}
\newtheorem{conj}[thm]{Conjecture}
\newtheorem{prop}[thm]{Proposition}

\theoremstyle{definition}
\newtheorem{defn}[thm]{Definition}

\theoremstyle{remark}
\newtheorem{rem}[thm]{Remark}

\numberwithin{equation}{section}
\numberwithin{thm}{section}

\begin{document}

\title{Extremal results on feedback arc sets in digraphs}

\author{Jacob Fox\thankssymb{1}}
\thanks{\thankssymb{1} Department of Mathematics, Stanford University, Stanford, CA 94305. Email: {\tt jacobfox@stanford.edu}. Research supported in part by a Packard Fellowship and by NSF Award DMS-1855635.}
\author{Zoe Himwich\thankssymb{2}}
\thanks{\thankssymb{2} Department of Mathematics, Columbia University, New York, NY. Email:
\href{mailto:himwich@math.columbia.edu} {\nolinkurl{himwich@math.columbia.edu}}.
}
\author{Nitya Mani\thankssymb{3}}
\thanks{\thankssymb{3} Department of Mathematics, Massachusetts Institute of Technology, Cambridge, MA. Email: \href{mailto:nmani@mit.edu} {\nolinkurl{nmani@mit.edu}}.
Research supported in part by a Hertz Fellowship and NSF GRFP}

\date{\today}

\maketitle
\begin{abstract}
A directed graph is \textit{oriented} if it can be obtained by orienting the edges of a simple, undirected graph. For an oriented graph $G$, let $\beta(G)$ denote the size of a \textit{minimum feedback arc set}, a smallest subset of edges whose deletion leaves an acyclic subgraph. A simple consequence of a result of Berger and Shor is that any oriented graph $G$ with $m$ edges satisfies $\beta(G) = m/2 - \Omega(m^{3/4})$.

We observe that if an oriented graph $G$ has a fixed forbidden subgraph $B$, the upper bound of $\beta(G) = m/2 - \Omega(m^{3/4})$ is best possible as a function of the number of edges if $B$ is not bipartite, but the exponent $3/4$ in the lower order term can be improved if $B$ is bipartite. We also show that for every rational number $r$ between $3/4$ and $1$, there is a finite collection of digraphs $\mcB$ such that every $\mcB$-free digraph $G$ with $m$ edges satisfies $\beta(G) = m/2 - \Omega(m^r)$, and this bound is best possible up to the implied constant factor. The proof uses a connection to Tur\'an numbers and a result of Bukh and Conlon. Both of our upper bounds come equipped with randomized linear-time algorithms that construct feedback arc sets achieving those bounds. Finally, we give a characterization of quasirandom directed graphs via minimum feedback arc sets.
\end{abstract}

\section{Introduction}
Throughout, let $G = (V, E)$ be a directed graph (abbreviated \textit{digraph}) such that for vertices $a, b \in V$, $(a, b) \in E$ is the edge directed $a \rightarrow b$. All digraphs in this paper are \textit{oriented} (i.e.~with no loops, parallel or antiparallel edges). A \textit{directed $r$-cycle} in $G$ consists of vertices $v_1,\ldots,v_r$ such that 
$(v_i,v_{i+1})$ is an edge for $i=1,\ldots,r$, where $v_{r+1}$ is taken to be $v_1$. A digraph is \textit{acyclic} if it has no directed cycles. 
In this article, we study extremal problems on the minimum number of edges we need to delete from a digraph to make the remaining digraph acyclic. 

\begin{defn} Given a digraph $G = (V, E)$, a \textit{feedback arc set} is a subset $S \subset E$ such that $G' = (V, E \backslash S)$ is acyclic. A \textit{minimum feedback arc set} is a feedback arc set of minimum size. Let $\beta(G)$ be the size of a minimum feedback arc set of $G$. 
\end{defn}

Understanding minimum feedback arc sets is of interest in mathematics and computer science. Feedback arc sets were originally studied by Slater~\cite{SLA61} in the context of \textit{tournaments}, oriented complete graphs, looking at inconsistent round-robin tournament outcomes. Problems related to minimum feedback arc sets arise in property testing, tearing chemical engineering, deadlock resolution, ranked voting, electronic circuits, and a number of other natural contexts (c.f. \cite{ELS93, LEI91}).
Unfortunately, computing the size of a minimum feedback arc set is \textsf{NP}-hard
even for tournaments (see \cite{AL06,CHAR07}). Consequently, there has been considerable work to show upper bounds on the size of minimum feedback arc sets as well as to develop efficient algorithms that can approximate $\beta(G)$ and/or construct provably small feedback arc sets (see~\cite{BS90,BS97,ENSS98,KST90,LR99}).

Observe that any oriented graph $G$ with $m$ edges satisfies the simple upper bound $\beta(G) \le m/2$. This bound follows by fixing any linear ordering $<$ of the vertices of $G$. We can partition the edge set of $G$ into two acyclic digraphs, one consisting of the ``forward'' edges, those $i \rightarrow j$ with $i < j$, and the other consisting of the ``backward'' edges, $k \rightarrow l$ with $k > l$. Deleting the smaller of these two sets yields an acyclic digraph with at least $m/2$ edges.

This straightforward upper bound $\beta(G) \le m/2$ illustrates a parallel between the study of $\beta(G)$ and the maximum cut of an undirected graph. When studying the extremal versions of both problems, we work to show how large an advantage over $m/2$ is possible, usually as a sublinear function of the number of edges of the associated oriented or undirected graph. 

Given the hardness of exactly computing a minimum feedback arc set, showing bounds on $\beta(G)$ for general digraphs and specific families of digraphs has been a topic of much interest since the 1960s.
Improving earlier work of Erd\H{o}s and Moon~\cite{EM65}, Spencer~\cite{SPE71,SPE80}, de la Vega~\cite{VEG83}, and Poljak, R\"odl, and Spencer~\cite{PRS88} showed sharp bounds for the size of a minimum feedback arc set in a tournament, observing that any tournament $T$ with $n$ vertices and $m={n \choose 2}$ edges satisfies $\beta(T) \le \frac{m}{2} - c m^{3/4}$ for a constant $c > 0$, and this bound is best possible up to the constant $c$ by considering a random tournament.

In a more general context, Berger and Shor~\cite{BS90,BS97} gave an algorithmic proof that for any digraph $G$ with $m$ edges and maximum degree $\Delta$ satisfies $\beta(G) \le \frac{m}{2} - c \frac{m}{\sqrt{\Delta}}$ for some constant $c > 0$. Their work also implies the following upper bound on the size of a minimum feedback arc set in a digraph and gives a randomized polynomial time algorithm to construct such a small feedback arc set. For a digraph, we let $d(v)$ denote the total degree of vertex $v$, defined as the sum of the indegree and outdegree of $v$.

\begin{thm}[\cite{BS97}]\label{t:betaub}
There is a constant $c>0$ such that for all digraphs $G = (V, E)$ with $m$ edges, 
$$\beta(G) \le \frac{m}{2}- c \sum_{v \in V} \sqrt{d(v)} \leq \frac{m}{2}- \frac{c}{4}  \cdot m^{3/4}.$$
\end{thm} 

Berger~\cite{BER97} used a fourth moment method to give a stronger upper bound that includes a contribution from the total discrepancy between the indegree and outdegree of the digraph. For completeness, we give a simple proof of Theorem~\ref{t:betaub} in the appendix. The second inequality is not directly in the earlier works but has a simple proof. 

As noted earlier (in work of~\cite{SPE71, SPE80, VEG83}), the upper bound on $\beta(G)$ given by Theorem~\ref{t:betaub} is best possible up to the constant factor in the lower order term as realized with high probability by a random tournament on $n$ vertices. The same proof method also implies that a random orientation of a balanced, complete bipartite graph with $m$ edges almost surely satisfies the same lower bound on $\beta(G)$.

\subsection*{Forbidden subgraphs} We might hope to improve upon the upper bound on the size of a minimum feedback arc set given by Theorem \ref{t:betaub} when we have certain further information about the digraph. For digraph $B$, let the \textit{underlying undirected graph} of $B$, denoted $\overline{B}$, be the graph where we replace all directed edges in $B$ with undirected edges.

We first take the case where $G$ is \textit{$B$-free} for some fixed forbidden digraph $B$. As we remarked earlier, a balanced complete bipartite graph with $m$ edges oriented uniformly at random with high probability has minimum feedback arc set of size at least $\frac{m}{2} - O(m^{3/4})$. Since bipartite graphs are $\overline{B}$-free for any non-bipartite $\overline{B}$,
it is primarily interesting to forbid some orientation of a bipartite graph. In this setting, we observe that if $B$ is a bipartite digraph whose underlying undirected graph is connected and which has edges in both directions between its parts, then there are $B$-free digraphs with quadratically many edges; indeed, the orientation of a complete bipartite graph with all edges directed from one part to another is $B$-free. This implies that in order to give an improved upper bound on $\beta(G)$, we cannot simply apply Theorem~\ref{t:betaub}. Nevertheless, in this case we can do substantially better, improving the exponent of $3/4$ in the lower order term.

\begin{thm}\label{t:betahfree}
For each bipartite digraph $B$, there are $c,\epsilon>0$ depending on $B$ such that the following holds. If digraph $G$ has $m$ edges and is $B$-free, then 
$$\beta(G) \le \frac{m}{2} - c m^{3/4 + \epsilon}.$$
\end{thm}
Sparse random digraphs (the edges and their orientation are random) can be used to show that Theorem~\ref{t:betahfree} is best possible apart from the exponential constant $\epsilon=\epsilon(B)$. The proof of Theorem \ref{t:betahfree} can be made algorithmic to find such a small feedback arc set with a randomized linear time algorithm for $B$ fixed. This is not immediately clear from the proof, as it uses a result in extremal graph theory related to Sidorenko's conjecture involving counts of subgraphs, but the use of sampling arguments allows one to accomplish this. Details are discussed in Remark~\ref{r:algbfree}. 

We also consider imposing more stringent conditions on a digraph $G$, by forbidding its underlying undirected graph $\overline{G}$ from containing any $B \in \mcB$, for some family $\mcB$ of undirected subgraphs. The \textit{extremal number} $\ext(n, \mcB)$ of $\mcB$ is the maximum number of edges a $\mcB$-free undirected graph on $n$ vertices can have. We obtain an improved upper bound on $\beta(G)$ under these conditions by relating $\beta(G)$ to $\ext(n, \mcB)$.

\begin{thm}\label{t:bfreeub}
Suppose $\mcB$ is a family of undirected graphs for which $\ext(n,\mcB)=O(n^{2-\eps(\mcB)})$. There is a constant $c=c(\mcB)>0$ such that for all digraphs $G = (V, E)$ with $m$ edges such that $\overline{G}$ is $\mcB$-free, we have $$\beta(G) \le \frac{m}{2} - c m^{\frac34 +  \frac{\eps(\mcB)}{4(2 - \eps(\mcB))}}.$$
\end{thm}

The proof of Theorem \ref{t:bfreeub} in fact shows that we obtain the bound from Theorem~\ref{t:betaub}, and hence by the Berger-Shor algorithm~\cite{BS97}, we can construct such a feedback arc set in randomized linear time. Further, the above bound is tight up to the constant coefficient of the directed surplus in the following sense.

\begin{thm}\label{t:bfreelb}
Let $\mcB$ be a collection of undirected graphs with $\ext(n, \mcB) = \Theta(n^{2 - \eps(\mcB)})$. 
There is a constant $C = C(\mcB) > 0$ such that for all positive integers $n$, there exists a digraph $G$ on $n$ vertices and $m = \ext(n, \mcB)$ edges where $\overline{G}$ is $\mcB$-free such that
$$\beta(G) \ge \frac{m}{2} -  C \cdot  m^{\frac34 + \frac{\eps(\mcB)}{4(2 - \eps(\mcB))}}.$$
\end{thm}

Notably, combined with the result of Bukh and Conlon~\cite{BC18}, Theorems \ref{t:bfreeub} and \ref{t:bfreelb} imply that the directed surplus achieves every rational exponent between $3/4$ and $1$, which we make more precise in Corollary~\ref{cor:hiteveryexp}.

\begin{rem}
The proof of Theorem \ref{t:bfreelb} uses a random orientation of any undirected graph that is $\mcB$-free with $n$ vertices and with  $m=\ext(n, \mcB)$ edges. The extra condition $\ext(n, \mcB) = \Theta(n^{2 - \eps(\mcB)})$ is not generally needed and we suspect such a random digraph $G$ likely satisfies $\beta(G)= m/2 - \Omega\left(\sqrt{mn}\right)$. The proof gives this as long as $\ext(n, \mcB)$ satisfies some mild conditions, and we only use the condition $\ext(n, \mcB) = \Theta(n^{2 - \eps(\mcB)})$ as it is a simple one to work with. 
\end{rem}

Forbidding cycles as subgraphs is of particular interest. A digraph is \textit{$r$-free} if it is free of all directed cycles of length at most $r$. Directed cycles in digraphs are far more poorly understood than cycles in graphs, with many fundamental problems remaining open. The archetypal example of the difficulty gap is the famous  Cacceta-H\"{a}ggkvist conjecture \cite{CAC78}, which states that all $r$-free digraphs on $n$ vertices have a vertex of outdegree less than $n/r$. This conjecture is still open even in the case $r = 3$, a stark contrast to the undirected version of this problem.

Previous work has largely focused on computing $\beta(G)$ in terms of $\gamma(G)$, the number
of non-adjacent, unordered pairs of vertices of $G$. Chudnovsky, Seymour and Sullivan~\cite{CHUD08} conjectured that if $G$ is triangle-free then $\beta(G) \le \gamma(G)/2$, which they showed in two special cases, showing the weaker bound $\beta(G) \le \gamma(G)$ for all digraphs $G$. This conjecture was generalized to \textit{$r$-free digraphs} by Sullivan~\cite{SUL08}. She conjectured that $\beta(G) \le 2 \frac{\gamma (G)}{(r+1)(r-2)}$ for all $r$-free digraphs $G$. This conjectured bound was proved up to a constant factor by Fox, Keevash, and Sudakov \cite{FOX10}, who showed that every $r$-free digraph $G$ satisfies $\beta(G) \le 800\gamma(G)/r^2$.%

Theorems~\ref{t:bfreeub} and~\ref{t:bfreelb} imply bounds in the simpler setting of digraphs whose underlying undirected graph has no $r$-cycles. Moving towards the general case, we construct a family of $r$-free digraphs $G$ on $n$ vertices with $$\beta(G) \ge \frac{n^2}{(r + 1)^2}.$$
This implies that if every $r$-free digraph $G$ on $n$ vertices satisfies $\beta(G) =O(1)$, then $r = \Omega(n)$. In the converse direction, we conjecture all $r$-free digraphs on $n$ vertices with $r > 2n/3$ are at most one edge away from being acyclic. If true, this would be tight, as for $n$ a multiple of $3$ there exists an $(r-1)$-free digraph $G$ on $n$ vertices with $r=2n/3$ and $\beta(G) = 2$.

\subsection*{Quasirandom directions and feedback arc sets} We are also able to use minimum feedback arc sets to characterize and study \textit{quasirandomness} in digraphs. Understanding when a graph can be deterministically shown to have ``random-like'' behavior has been a problem that has been studied for a long time, first by Thomason~\cite{THO87} and systematically by Chung, Graham, and Wilson~\cite{CGW89}.

A precise definition of quasirandomness in undirected graphs along with several equivalent characterizations were given by Chung, Graham, and Wilson~\cite{CGW89}. These ideas were generalized by Chung and Graham to notions of quasirandomness for tournaments~\cite{FAN91} and subsequently directed and oriented graphs in later work~\cite{GRI13, AM11}. Notably, Griffiths~\cite{GRI13} studied whether a digraph is oriented quasirandomly with respect to its underlying undirected graph, the notion of quasirandomness we develop further here.

The definition of \textit{quasirandom direction} we use will rely on the following notion of the maximum discrepancy in edges between two subsets of a digraph.

\begin{defn} \label{d:discrepancy}
Given a digraph $G = (V, E)$, for two subsets $A,B\subset V$ we define the \textit{directional discrepancy} of $G$ as $$\tau(G):= \max_{A,B\subset V}(e(A,B)-e(B,A)),$$
where $e(A, B) = |\{(a, b) \in E : a \in A, b \in B\}|$.
\end{defn}

Throughout our discussion of quasirandomness in this article, we examine a set of properties that a dense digraph $G = (V, E)$ on $n$ vertices and $m = \Omega(n^2)$ edges might satisfy when describing notions of quasirandomness.  

\begin{defn}
A digraph $G$ with underlying undirected graph $\overline{G}$ and $m$ edges has \textit{quasirandom direction} with respect to $\overline{G}$ if $\tau(G) = o(m)$.
\end{defn}

We can understand whether a dense digraph is quasirandom by looking at the the size of its minimum feedback arc set.
\begin{thm}\label{t:qrandom}
For a digraph $G = (V, E)$ on $n$ vertices and $m = \Omega(n^2)$ edges with underlying undirected graph $\overline{G}$, $G$ has quasirandom direction with respect to $\overline{G}$ if and only if $\beta(G) = m/2 - o(m)$.
\end{thm}
For completeness, we give an expanded characterization of quasirandom direction in Section~\ref{s:qrandom}, including many properties which are directly found in~\cite{GRI13,AM11} or arise as straightforward consequences of the techniques from earlier related work on quasirandomness by e.g.~\cite{FAN91}.

\section{Preliminaries}\label{s:boundbeta}
We first introduce some notation used throughout the article.
Given digraph $G = (V, E)$ and subsets $A,B\subset V$, let $$e(A, B) = |\{e = (a, b) \in E \mid a \in A,\,b\in B\}|.$$
For a vertex $v \in V$, let $d^+(v)$ be the \textit{indegree} of a vertex $v \in V$, i.e. 
$$d^+(v) = | \{ w \in V \mid (w, v) \in E \}|,$$
and similarly let $d^-(v)$ denote the \textit{outdegree} of $v$ (recall that $d(v) = d^+(v) + d^-(v)$ is the total degree of $v$). 
When computing properties of feedback arc sets for digraphs, it is often helpful to fix an ordering of the vertices. 

\begin{defn}
Given a digraph $G = (V, E)$ on $n$ vertices, an \textit{ordering} of the vertices is a bijection $\rho: V \rightarrow [n]$. Given an ordering of $V$, an edge $e = (v, w) \in E$ is a \textit{forward edge} if $\rho(v) < \rho(w)$, else it is a \textit{backwards edge}.
\end{defn}

Younger~\cite{YOU63} showed that a digraph is acyclic if and only if there exists an ordering of its vertex set for which there are no backward edges. 

In this work, we study the following extremal question: given a digraph with $m$ edges and some constraints on the digraph, how few edges can we remove to guarantee that the resulting digraph is acyclic? 

We observe that there is a close relationship between the size of a minimum feedback arc set and the directional discrepancy of a digraph (recall Definition~\ref{d:discrepancy}). We will be interested in two additional variants of the directional discrepancy.

\begin{defn} \label{d:discrepancy2}
Given a digraph $G = (V, E)$ let
$$\tau^{*}(G):= \max_{A,B\subset V, A\cap B=\emptyset}(e(A,B)- e(B,A)).$$ 
We also define the \textit{maximum edge difference for partitions}:
$$\tau_{\sqcup}(G) := \max_{A \sqcup B = V} (e(A,B)- e(B,A)).$$
\end{defn}

Note that $\tau^*(G)$ is not necessarily achieved by a partition and thus we have the inequality $\tau_{\sqcup}(G) \le \tau^*(G)$, but not necessarily the reverse inequality as shown by a cyclic triangle. However, $\tau(G)$ and $\tau^*(G)$ are within a constant factor of each other.

\begin{lemma} \label{l:tautau}
For any digraph $G$, we have $\tau^*(G) \le \tau(G) \le 3 \tau^*(G)$.  
\end{lemma}
\begin{proof}
By definition, $\tau(G) \ge \tau^*(G)$. It remains to show the second inequality.  Choose $A, B \subset V$ so that $\tau(G) = e(A, B) - e(B, A).$ Let $C = A \cap B$.
We have that
\begin{align*}
\tau(G) &= e(A, B) - e(B, A) \\
&= (e(A  \backslash C, B  \backslash C) + e(C, B  \backslash C) + e(A  \backslash C, C) + e(C, C) ) \\
&\qquad -  (e(B  \backslash C, A  \backslash C) + e(C, A  \backslash C) + e(B  \backslash C, C) + e(C, C) ) \\
&= (e(A  \backslash C, B  \backslash C) - e(B  \backslash C, A  \backslash C) ) + (e(C, B  \backslash C) - e(B  \backslash C, C) ) \\
&\qquad + (e(A  \backslash C, C) ) -  (e(C, A  \backslash C)) \\
&\le 3 \tau^*(G).
\end{align*}
\end{proof}
In light of the simple upper bound $\beta(G) \leq e(G)/2$, the following definition is natural. 

\begin{defn} The \textit{directed surplus} of a digraph $G$ on $m$ edges is $$\pi(G):=\frac{m}{2}-\beta(G).$$
\end{defn}

The directed surplus is half of the maximum difference between the number of forward and backward edges over all orderings of the vertices of $G$. A simple lemma shows that $\pi(G)$ is monotone on induced subgraphs.

\begin{lemma}\label{l:monotone}
For any digraph $G = (V, E)$ and vertex subset $U \subset V$, we have 
$$\pi(G) \ge \pi(G[U]),$$
where $G[U]$ is the induced subdigraph of $G$ on $U$.
\end{lemma}
\begin{proof}
Take an optimal ordering $\rho$ of $U$ in digraph $G[U]$. We extend $\rho$ to an ordering of $V(G)$ by uniformly at random adding each unseen vertex either to the front or back of the ordering. In expectation, half of the newly added edges are backwards and thus, by choosing some ordering with at least as many new forwards edges as backwards edges, the surplus of this ordering on $V(G)$ is at least as large as the one on $G[U]$. Hence, $\pi(G) \ge \pi(G[U])$.
\end{proof}

We can upper bound the directional discrepancy $\tau(G)$ by a function of the directed surplus of digraph $G$.

\begin{lemma}\label{t:betatau}
For any digraph $G = (V, E)$,
$$\frac{1}{6} \tau(G) \le \pi(G).$$
\end{lemma}
\begin{proof}
Via Lemma~\ref{l:tautau}, it suffices to show that $\frac12 \tau^*(G) \le \pi(G)$.
Pick disjoint $A$ and $B$ so that $\tau^*(G) = e(A, B) - e(B, A)$. We take an ordering of each of $A$ and $B$ to maximize the number of forwards edges within $A$ and $B$. This induces an ordering $\rho$ on $A\sqcup B$, so that for all $a \in A$ and all $b \in B$, $\rho(a) < \rho(b)$. In this ordering, the difference between the number of forwards and backwards edges is at least $\tau^*(G)$, as within $A, B$ we have at least as many forwards as backwards edges. 

We extend $\rho$ to an ordering on $V$ greedily. One by one, we assign each $v \in V \backslash (A \cup B)$ to be either the first or last vertex in the ordering so that the difference between forwards and backwards edges in $\rho$ is non-decreasing. In doing this, the final ordering $\rho$ has at least $\tau^*(G)$ more forwards edges than backwards edges, which implies the desired bound.
\end{proof}

\section{Digraphs with a forbidden bipartite subgraph} 
The goal of this section is to prove Theorem \ref{t:betahfree}, giving a better upper bound on $\beta(G)$ for $B$-free digraphs $G$ where $\overline{B}$ is bipartite. 
We show that any two orientations of a fixed undirected subgraph $\overline{B}$ have similar numbers of copies in digraph $G$ if $\tau(G)$ is relatively small (see Corollary~\ref{c:smalltau}). An approximate version of Sidorenko's conjecture implies that $\overline{G}$ must contain many copies of $\overline{B}$ (as observed in Lemma~\ref{l:bdtau}). Thus for there to be no copy of digraph $B$, either $\tau(G)$ must be large (in which case $\beta(G)$ is relatively small) or $G$ cannot have too many edges. We handle the cases when digraph $G$ is relatively sparse using Theorem~\ref{t:betaub}.

\begin{defn}
For digraphs $G$ and $B$, let $N_L(B,G)$ be the number of labeled copies of $B$ in $G$. Similarly, if $H$ and $F$ are undirected graphs, let $N_L(F,H)$ be the number of labeled copies of $F$ in $H$.
\end{defn}

\begin{prop} \label{p:countflip1}
Let $B$ be a digraph and let $B_1$ be the digraph formed by reversing a single edge $e = (u, v)$ of $B$. Let $D$ be the induced subgraph of $B$ on $V(B) \backslash \{u, v\}$. Then for any digraph $G$,
$$| N_L(B, G) - N_L(B_1, G)| \le \tau(G) \cdot N_L(D, G).$$
\end{prop}
\begin{proof}
Fix a copy $D_0$ of $D$ in $G$. Let $S_u$ be the set of $w \in V(G)$ so that the induced graph on $V(D_0) \cup \{w\}$ contains a copy of $B[V(B) \backslash \{v\}]$ as a labeled subgraph (where $u$ is the image of $w$).
Similarly, let $S_v$ be the set of $w \in V(G)$ so that the induced graph on $V(D_0) \cup \{w\}$ is isomorphic to $B[V(B) \backslash \{u\}]$. 
The number of ways to extend $D_0$ to a labeled subgraph $B$ in $G$ is $e(S_u,S_v)$. By definition of $B_1$, the number of ways to extend $D_0$ to a copy of $B_1$ in $G$ is exactly $e(S_v, S_u)$. Thus, considering all possible copies of $D$ in $G$, by the triangle inequality and definition of $\tau(G)$, we obtain the desired bound as follows: 
$$|N_L(B, G) - N_L(B_1, G) | \le N_L(D, G) \cdot \max_{D_0 \subset G}|e(S_u, S_v) - e(S_v, S_u)| \le N_L(D, G) \cdot \tau(G).$$
\end{proof}

\begin{cor}\label{c:smalltau}
If $B, B'$ are digraphs with the same underlying undirected graph $\overline{B}$ on $k$ vertices and $G$ is a digraph on $n$ vertices, then
$$| N_L(B, G) - N_L(B', G)| \le e(B) \cdot n^{k-2} \cdot \tau(G).$$
\end{cor}
\begin{proof}
By flipping the direction of at most $e(B)$ edges, we obtain a copy of $B'$ from a copy of $B$ (since $B, B'$ have the same underlying undirected graph). Thus, by applying Proposition~\ref{p:countflip1} $e(B)$ times and the triangle inequality, we find that 
$$|N_L(B, G) - N_L(B', G)| \le e(B) \cdot \tau(G) \cdot \max_{D \in \mcR} N_L(D, G),$$
where $\mcR$ is the set of labeled subgraphs obtained by removing exactly one edge of any orientation of $\overline{B}$. Since $N_L(D, G) \le n^{k-2}$ for all such digraphs $D$, we obtain the desired bound:
$$|N_L(B, G) - N_L(B', G)| \le n^{k-2} \cdot e(B) \cdot \tau(G).$$
\end{proof}

\begin{rem}
The bound $N_L(D,G) \leq n^{k-2}$ in the above proof, and hence the factor $n^{k-2}$ in the bound in Corollary~\ref{c:smalltau}, can be improved for many $D$ if $G$ is sparse using a result of Alon~\cite{AL81}, which would give an improvement on the $\eps(B)$ obtained in the proof of Theorem~\ref{t:betahfree}. 
\end{rem}

The above result has the useful consequence that if digraph $G$ is $B$-free and $\overline{G}$ has many copies of $\overline{B}$, then $\tau(G)$ is small. If $\tau(G)$ is large, since all orientations of $\overline{B}$ are roughly equally represented in $G$, $G$ will contain a copy of $B$. We will use this interplay between the size of $\tau(G)$ and Theorem~\ref{t:betaub} to improve our upper bound on $\beta(G)$ when $G$ is $B$-free.
We split our analysis into two cases; we first show a relationship between $m$ and $\tau(G)$ by leveraging Sidorenko's conjecture.

\begin{defn}
A \textit{homomorphism} from an undirected graph $H$ to an undirected graph $G$ is a mapping $\rho: V(H) \rightarrow V(G)$ such that, for each $(u, v) \in E(H)$, $(\rho(u), \rho(v)) \in E(G)$ (i.e. $\rho$ maps edges to edges). Let $h_H(G)$ denote the number of homomorphisms from $H$ to $G$, and let $t_H(G) := h_H(G)/|V(G)|^{|V(H)|}$ be the normalized homomorphism count, i.e. the fraction of vertex maps which are homomorphisms.
\end{defn}

The \textit{edge density} of a graph $G$ with $n$ vertices and $m$ edges is $t_{K_2}(G) = 2m/n^2$.
A famous conjecture of Sidorenko \cite{Si93} and Erd\H{o}s-Simonovits \cite{ES83} states that for fixed $0<p<1$, the random graph $G(n,p)$ asymptotically minimizes (over all graphs on $n$ vertices with edge density $p$) the number of copies of a fixed bipartite graph $H$. The following is a more precise version in terms of homomorphism densities.

\begin{conj}[Sidorenko's Conjecture]
For every undirected bipartite graph $H$ with $e(H)$ edges and every undirected graph $G$,
$$t_H(G) \ge t_{K_2}(G)^{e(H)}.$$
\end{conj}

Conlon, Fox, and Sudakov~\cite{CON10} showed that Sidorenko’s conjecture holds for every undirected bipartite graph $H$ which has a vertex complete to the other part. The {\it width} $w(H)$ of an undirected bipartite graph $H$ is the minimum number of edges that must be added to $H$ to obtain a bipartite graph where a vertex of $H$ is complete to the other part (so at most $v(H)/2$). For an undirected bipartite graph $H$, let the \textit{exponent} of $H$ be $t(H)=e(H)+w(H)$.
Their work implies the following quantitative bound for all bipartite $H$.
\begin{lemma}[Corollary 1.1~\cite{CON10}]\label{cor:manycopies}
Fix undirected bipartite graph $H = (U_1 \sqcup U_2, F)$ with $k$ vertices and exponent $t$. The number of labeled copies of $H$ in any undirected graph $G = (V, E)$ with $n$ vertices and edge density $p \ge (k^2/n)^{1/t}$ is at least $\frac12 n^k p^{t}$.
\end{lemma}
The density condition implies that the number of labeled copies of $H$ in $G$ dominates the number of non-injective maps (and hence homomorphisms) of $H$ into $G$.
It is noteworthy that several further cases of Sidorenko's conjecture have been solved~\cite{CL17,CL21,LI11,CON18,H10,SZE14,KIM16}.

In the direction of proving Theorem \ref{t:betahfree}, we first observe that $B$-free digraphs $G$ (where $B$ has bipartite underlying undirected graph) have large $\tau(G)$.

\begin{lemma}~\label{l:bdtau} 
Let $B$ be a directed bipartite graph with $k$ vertices such that $\overline{B}$ has exponent $t$. If a digraph $G = (V, E)$ on $n$ vertices and $m \ge \frac12 k^{2/t} n^{2-1/t}$ edges is $B$-free, then we have $\tau(G) \ge c \cdot \frac{m^t}{n^{2t - 2}}$ for some $c=c(B)>0$. 
\end{lemma}
\begin{proof}
The assumption giving a lower bound on $m$ is equivalent to the edge density condition in Lemma~\ref{cor:manycopies}. Hence, Lemma~\ref{cor:manycopies} implies that, letting $p=2m/n$ denote the edge density of $G$,  
$$N_L(\overline{B}, \overline{G}) \ge \frac12 n^k p^{t}.$$
Suppose $B’$ is a digraph with the same underlying undirected graph as $B$. By Corollary~\ref{c:smalltau}, we have 
$$| N_L(B, G) - N_L(B’, G)| \le e(B) n^{k-2} \tau(G).$$  
Since $G$ is $B$-free and we can orient the edges of $\overline{B}$ in at most $2^{e(B)}$ ways, $\overline{G}$ can have at most 
$(2^{e(B)}-1)e(B)n^{k-2}\tau(G)$ copies of $\overline{B}$. Therefore,  
$$(2^{e(B)}-1)e(B)n^{k-2}\tau(G) \ge N_L(\overline{B}, \overline{G}) \ge \frac12 n^k p^{t}.$$
This implies the desired bound on $\tau(G)$:
$$\tau(G) \ge \frac{p^t}{2 \cdot (2^{e(B)}-1) \cdot e(B)} n^2 \ge \frac{p^{t}}{2^{e(B)+1} \cdot e(B)} \cdot n^2 \ge \frac{1}{2e(B)} \cdot \frac{m^{t}}{n^{2t-2}}.$$
\end{proof}

\begin{proof}[Proof of Theorem~\ref{t:betahfree}]
Let $t=t(B)$ and let digraph $G$ have $n$ vertices and $m$ edges. Let $\alpha=\frac{1}{2t-1}$ and $\eps=\frac{1}{16t-12}$. The proof splits into two cases. 

In the first case, $m >n^{2-\alpha}$. We apply Lemma~\ref{l:bdtau}, which in conjunction with Lemma~\ref{t:betatau}, that $\pi(G) \ge \frac16 \tau(G)$, implies that for some $c_1 = c_1(B) > 0$,
$$\pi(G) \geq \frac16 \tau(G) \ge c_1\frac{m^t}{n^{2t-2}} > c_1m^{t-(2t-2)/(2-\alpha)}  = c_1 m^{\frac34 + \eps},$$
and hence $\beta(G)=\frac{m}{2}-\pi(G) > \frac{m}{2}-c_1 m^{\frac34 + \eps}$.

In the second case, we have $m \le n^{2-\alpha}$. Let $U \subset V(G)$ be the vertices of degree at least $m^{1/2 - 2\eps}$. Note that $|U| \leq 2m/m^{1/2 - 2\eps}=2m^{1/2 + 2\epsilon}$. 

If $e(G[U]) \ge m/2$, then we apply Lemma~\ref{t:betatau} and then Lemma~\ref{l:bdtau} analogously to above but to $G[U]$, which gives that for some $c_2 = c_2(B) > 0$,
$$\pi(G[U]) \ge  \frac16 \tau(G[U]) \ge \frac{c}{6} \cdot  \frac{e(G[U])^t}{|U|^{2t - 2}} \ge \frac{c}{6} \cdot  \frac{\left(m/2\right)^t}{\left(2m^{1/2 + 2\epsilon}\right)^{2t - 2}}  \ge c_2 m^{\frac34 + \eps}.$$
By Lemma~\ref{l:monotone}, we then find $\beta(G) \le \beta(G[U]) \le \frac{m}{2} - c_2 m^{\frac34 + \eps}$.

It remains to consider the case $e(G[U]) < m/2$, that is, where more than half of the edges of $G$ are incident to vertices of degree at most $m^{1/2 - 2\eps}$. Restricting our attention to low degree vertices, we see that 
$$\sum_{v \in V} \sqrt{d(v)} \ge \sum_{v : d(v) \le m^{1/2 - 2\eps}} \frac{d(v)}{\sqrt{d(v)}} \ge \frac{m/2}{\sqrt{m^{1/2 - 2\eps}}} = \frac12 m^{3/4 + \eps}.$$
This implies the desired bound by Theorem~\ref{t:betaub}:
$$\beta(G) \le\frac{m}{2} - c \sum_{v \in V}\sqrt{d(v)} \le  \frac{m}{2} - \frac{c}{2} m^{\frac{3}{4} + \eps}.$$

The result then follows with $\eps(B) =\eps = \frac{1}{16t-12}$ and with the constant factor chosen to be $\min(c/2, c_1, c_2)$.
\end{proof}

\begin{rem}\label{r:algbfree}
We observe that the above proof of Theorem~\ref{t:betahfree} can be adapted to get a randomized linear time algorithm for constructing a small feedback arc set; we sketch this algorithm below. First, if $m \leq n^{2-\alpha}$, the proof above implies that the Berger-Shor algorithm~\cite{BS97} achieves the desired small feedback arc set.

Thus, it remains to study denser graphs $G= (V, E)$ with $m > n^{2-\alpha}$. Recall $p=2m/n^2>n^{-\alpha}$. We will algorithmically find subsets $S, T \subset V$ such that $\tau(G) \ge e(S, T) - e(T, S) \ge c n^2 p^t$ for some $c = c(B) > 0$; we will use this pair of subsets to construct an ordering that yields a small feedback arc set.

By Lemma~\ref{cor:manycopies}, $\overline{G}$ contains at least $\frac12 n^{v(B)} p^t$ copies of $\overline{B}$. By averaging (as in the proof of Lemma~\ref{l:bdtau}), this implies that there are two orientations $B', B''$  of $\overline{B}$ that differ only in the orientation of a single edge $(u, v) \in E(B)$ such that for some $c' = c'(B) > 0$,
$$|N_L(B', G) - N_L(B'', G)| \ge c' n^{v(B)} p^t.$$
Let $D$ be the induced subgraph of $B'$ formed by deleting vertices $u$ and $v$. Since each copy of $D$ can extend in at most $n^2$ ways to an oriented copy of $\overline{B}$ in $G$, and the difference in the number of copies of $B'$ and $B''$ in $G$ is at least $c'n^{v(B)} p^t$, by Markov's inequality there must be at least $\frac{1}{2}c'n^{v(B)-2} p^t$ copies of $D$ in $G$ such that the difference in the number of copies of $B'$ and the number of copies of $B''$ that extend from this copy of $D$ is at least $\frac{1}{2}c'n^2p^t$.  

For a given copy of $D$ in $G$, let $S_u \subset V$ be the set of vertices that extend $D$ to a copy of  $B' \backslash \{v\}$, and let $S_v \subset V$ be the set of vertices that extend $D$ to a copy of $B' \backslash \{u\}$. This copy of $D$ together with the vertices of any edge $(u',v')$ with $u' \in S_u$ and $v' \in S_v$ form an extension of $D$ to $B'$. Similarly, any edge $(v',u')$ with $u' \in S_u$ and $v' \in S_v$ form an extension of $D$ to $B''$. Thus, the difference in the number of copies of $B'$ and the number of copies of $B''$ that extend from this copy is $|e(S_u, S_v) - e(S_v, S_u)|$.

Consider a fixed copy of $D$. It follows from the Chernoff bound that if we randomly sample $s=100\epsilon^{-1} \log n$ uniformly random pairs of vertices, with probability at least $1-n^{-2}$, the difference between the number of pairs that are are edges in $E(S_u,S_v)$ and the number of those that are edges in $E(S_v,S_u)$ differs from the expected value by at most $\epsilon s/2$. Thus, we can randomly sample a small number of pairs of vertices to detect if the absolute difference in the number of extensions of $D$ to $B'$ and to $B''$ is large.

Thus, by sampling subsets of $V$ of size $v(B) - 2$, after picking $O_B(n^{t \alpha})$ such samples, with probability at least $0.99$, we will find some subset that forms a copy of $D$ with the desired difference in the number of extensions that form copies of $B'$ and $B''$. By considering the associated sets $S_u, S_v$ to this copy of $D$, we will thus find subsets $S_u, S_v \subset V$ such that 
$$|e(S_u, S_v) - e(S_v, S_u)| \ge c(B) n^2 p^t.$$ In each iteration, we sample $O_B(n^{t\alpha}\log n)$ pairs to determine with high probability if $|e(S_u, S_v) - e(S_v, S_u)|$ is large, giving a run time of $$O_B(n^{2t \alpha}\log n)=O_B(n^{2t/(2t-1)}\log n)=o_B(m)$$ 
as $\alpha=1/(2t-1)$ and $m>n^{2-\alpha}$. 

Then, we can take an ordering $\rho$ on $G[S_u \cup S_v]$ so that $\rho(x) < \rho(y)$ for all $x \in S_u, y \in S_v$ and so that less than half the edges in $G[S_u], G[S_v]$ are backwards with respect to $\rho$. This latter condition can be obtained by taking an arbitrary relative ordering of $\rho$ on $S_u$ and either keeping or reversing that orientation to make at most half the edges of $G[S_u]$ be backwards (and analogously for $G[S_v]$ independently). Applying the algorithm in Lemma~\ref{l:monotone}, we find an associated small feedback arc set by extending $\rho$ greedily to an ordering on $V$ and deleting backwards edges with respect to $\rho$.
\end{rem}

\section{Forbidding families of undirected subgraphs}\label{s:boundbfree}
In this section, we study the size of a minimum feedback arc set of digraphs $G$ whose underlying undirected graph is $\mcB$-free for some finite collection of undirected subgraphs $\mcB$. We first present an upper bound on $\beta(G)$. This gives a lower bound on the directed surplus $\pi(G)$ as a function of the number of edges of $G$ which in many instances is best possible up to a constant factor. 

\subsection{Upper bound}
We show that for digraphs $G$ whose underlying undirected graph is $\mcB$-free, we can obtain an improved lower bound on the directed surplus of $G$ in terms of the order of $\ext(n, \mcB)$. We first observe a simple lower bound on $\sum_v \sqrt{d(v)}$ for digraphs where small induced subgraphs do not contain too many edges.

\begin{lemma}\label{l:manylowdeg1}
Suppose $G$ is a digraph with $m$ edges such that every induced subgraph of $G$ on $k$ vertices contains at most $m/2$ edges. Then 
$\sum_{v \in V(G)} \sqrt{d(v)} \geq \frac{1}{4}\sqrt{mk}$. 
\end{lemma}
\begin{proof}
Let $U$ be the $k$ vertices of $G$ of highest degrees and $W=V \setminus U$. Then each vertex in $W$ has degree at most $2m/k$. 
As $G[U]$ contains at most $m/2$ edges, and every edge not in $G[U]$ is incident to at least one vertex in $W$, we obtain $\sum_{v \in W} d(v) \geq m/2$.
By concavity of the function $f(x)=\sqrt{x}$, we have 
$$\sum_{v \in V} \sqrt{d(v)} \geq \sum_{v \in W} \sqrt{d(v)} \geq \left(\frac{m/2}{2m/k}\right) \cdot \sqrt{2m/k} \geq \frac{1}{4}\sqrt{mk}.$$
\end{proof}

Let $f(m,\mcB)$ denote the minimum number $n$ of vertices a $\mcB$-free graph with at least $m$ edges can have. Note that this is essentially the inverse function of $\ext(n,\mcB)$. Combining Theorem~\ref{t:betaub} and Lemma~\ref{l:manylowdeg1}, we have the following immediate corollary. 

\begin{cor}\label{cor:manylowdegext}
There is a constant $c>0$ such that the following holds. If a digraph $G$ has $m$ edges and the undirected graph $\overline{G}$ is $\mcB$-free, then $\beta(G) \leq m/2 - c\sqrt{mk}$ where $k=f(m/2,\mcB)$.
\end{cor}

Suppose $\mcB$ is such that we have $\ext(n,\mcB)=O(n^{2-\epsilon(B)})$ for all $n$. Then $f(m,\mcB)=\Omega(m^{1/(2-\epsilon(B))})$ for all $m$. Together with Corollary~\ref{cor:manylowdegext}, this directly implies Theorem~\ref{t:bfreeub}.

\subsection{Lower bound}
We construct a family of digraphs $G$ whose underlying undirected graphs are $\mcB$-free that meet the upper bound of Theorem~\ref{t:bfreeub} up to a constant factor of the directed surplus. These digraphs are obtained from an extremal $\mcB$-free graph (so with $n$ vertices and $\ext(n, \mcB)$ edges) by randomly orienting the edges. 

\begin{lemma}\label{lem:chernorient} 
Let $H$ be an undirected graph on $n$ vertices. There is a digraph $G$ obtained from $H$ by orienting its edges such that for every $s \in [n]$ and for all disjoint $A, B \subset V$ with $|A| = |B| = s$, we have that $$\left|e_G(A, B) - \frac12 \overline{e}(A, B) \right| \le 3\sqrt{\overline{e}(A, B)} \sqrt{s \log(en/s)},$$
where $\overline{e}(A, B)$ is the number of edges of $H$ with one endpoint in each of $A, B$.
\end{lemma}
\begin{proof}
Consider a random digraph $G$ obtained by uniformly at random orienting the edges of $H$. 
Consider disjoint vertex subsets $A,B \subset V$ with $|A|=|B|=s$. By the Chernoff bound (see Corollary A.1.7 in \cite{AlSp}), we have
$$\P\left(\left|e_G(A, B) - \frac12 \overline{e}(A, B) \right| > \delta \overline{e}(A, B) \right) < 2 \exp\left(- \frac{2}{3}\delta^2 \overline{e}(A, B)\right).$$
There are at most ${n \choose s}^2 \le (en/s)^{2s}$ ways to pick disjoint $A, B \subset V$ with $|A| = |B| = s$. Thus, by letting $\delta := \frac{3}{\sqrt{\overline{e}(A, B)}} \sqrt{s \log (en/s)}$ and taking a union bound, we observe that for a given positive integer $s$, the probability that \textit{any} $A, B \subset V$ with $|A| = |B| = s$ has $\left|e_G(A, B) - \frac12 \overline{e}(A, B) \right| > \delta \overline{e}(A, B)$ is at most 
$$\left(\frac{en}{s}\right)^{2s} \cdot 2 \exp(- 6s \log(en/s)) = 2 \left(\frac{en}{s} \right)^{-4s} \le \frac{1}{n}.$$ 
Taking a union over the at most $n/2$ choices of $s$, with probability at least $\frac12$, we have that for every positive integer $s$ and disjoint $A, B \subset V$ with $|A| = |B| = s$, the following inequality holds:
$$\left|e_G(A, B) - \frac12 \overline{e}(A, B) \right| \le \delta \overline{e}(A, B) = 3\sqrt{\overline{e}(A, B)} \sqrt{s \log(en/s)}.$$
\end{proof}

The above lemma allows us to construct a family of digraphs that meets the bound on $\beta(G)$ of Theorem~\ref{t:bfreeub} up to a constant factor of the directed surplus by considering a random orientation of an extremal $\mcB$-free undirected graph.

\begin{proof}[Proof of Theorem~\ref{t:bfreelb}]
Let $\mcB$ be a finite collection of undirected graphs with $\ext(n, \mcB) = \Theta(n^{2 - \eps(\mcB)})$. 
First, if $\mcB$ does not contain any forests, then $\eps(\mcB) < 1$, since for any integer $g > 2$ and any sufficiently large $n$, there exists a graph on $n$ vertices with girth $g$ and at least $n^{1+\eps_g}$ edges with $\eps_g > 0$ (c.f.~Corollary 2.28 in \cite{FS13}).
Else if $\eps(\mcB) = 1$, then with $C=1/2$ the inequality we are trying to show reduces to the trivial $\beta(G) \ge 0$.
Hereafter, we assume that $\eps(\mcB) < 1$. 

As long  as $n$ is at least the number of vertices in every graph in $\mcB$, a digraph obtained by adding isolated vertices to a $\mcB$-free digraph on $n$ vertices is also $\mcB$-free. Further, isolated vertices do not affect the size of a digraph's minimum feedback arc set. Note that for any positive integer $n$, there is a power of $2$ that is at least $n$ and less than $2n$. Thus by adding dummy isolated vertices, we can assume for convenience that $n$ is a power of $2$.

For any $n$, fix some extremal $\mcB$-free undirected graph with $n$ vertices and $m=\ext(n, \mcB)$ edges and let digraph $G$ be obtained by choosing an orientation of the edges that satisfies Lemma~\ref{lem:chernorient}.
Fix an arbitrary ordering $\rho$ on $V(G) = \{1, 2, \ldots, n\}$. 
A dyadic interval consists of the $2^i$ integers $(k2^i,(k+1)2^i]$ between two consecutive multiples of a power of $2$. Partition the set of pairs of distinct vertices so that each pair of vertices has one vertex in each of two consecutive dyadic intervals $A=\left(k2^i,(k+1)2^i\right]$ and $B=\left((k+1)2^i,(k+2)2^{i}\right]$ with $k$ even. We apply the bound above for $i=0,1,\ldots,\log_2 (n/2)$ between each such pair $A,B$ of consecutive dyadic intervals of size $2^i$. There are $n2^{-i-1}$ such pairs for a given $i$; we call this family of pairs $F_i$. Let $m_i$ be the number of edges between consecutive parts of size $2^i$ in the dyadic partition of $V(G)$ with respect to $\rho$ described above. 
Hence, $m = \sum_{i = 0}^{\log_2(n/2)} m_i$. 

We then use the triangle inequality to find that the number of forward edges with respect to ordering $\rho$ differs in absolute value from $m/2$ by at most 
\begin{align*}
\pi(G)&\le \sum_{i = 0}^{\log_2(n/2)} \sum_{\substack{(A, B)  \in F_i}} 3\sqrt{\overline{e}(A, B)} \sqrt{2^i \log(en/2^i)}. \\
\end{align*}
Since $G$ has $\mcB$-free underlying undirected graph, so does $G[A \cup B]$ and thus $\overline{e}(A, B) \le C(2^{i+1})^{2 - \eps(\mcB)}$ for an appropriate constant 
$C=C(\mcB)$.
Then, 
\begin{align*}
\pi(G)&\lesssim \sum_{i = 0}^{\log_2(n/2)} \sum_{\substack{(A, B)  \in F_i}} \sqrt{\overline{e}(A, B)} \sqrt{2^i \log(en/2^i)}. \\
&\lesssim \sum_{i = 0}^{\log_2(n/2)} n 2^{-i-1} \cdot  \sqrt{2^{i(3 - \eps(\mcB))} \cdot \log(en/2^i)}. \\
&\lesssim n \sum_{i = 0}^{ \log_2(n/2)} 2^{\left(1-\eps(\mcB)\right)i/2} \sqrt{ \log(en/2^i)}. \\
&\lesssim n \cdot 2^{(1 - \eps(\mcB)) \cdot \log (n/2) / 2 } \\
&\lesssim n^{(3 - \eps(\mcB))/2} \\
&\lesssim m^{\frac34 + \frac{\eps(\mcB)}{4(2 - \eps(\mcB))}}.
\end{align*}
This implies the desired lower bound on $\beta(G)$.
\end{proof}

A result of Bukh and Conlon~\cite{BC18} shows that there are finite families of undirected graphs with extremal number achieving every rational exponent between $1$ and $2$. 

\begin{thm}[Theorem 1.1~\cite{BC18}]
For every rational $1 \le r \le 2$, there exists a finite family of graphs $\mcB_r$ with $\ext(n, \mcB_r) = \Theta(n^r)$.
\end{thm}

Combining the above result with Theorems~\ref{t:bfreelb} and~\ref{t:bfreeub} shows that the directed surplus achieves every rational exponent between $3/4$ and $1$. The exponent $3/4+q$ is achieved by taking $r=8q/(4q+1)$ in the previous theorem. 
\begin{cor}\label{cor:hiteveryexp}
For every rational number $0 \leq q \leq 1/4$, there exists a finite family $\mcB$ of graphs and constants $c_1,c_2>0$ such that the following holds. For every digraph $G$ with $m$ edges whose underlying undirected graph  $\overline{G}$ is $\mcB$-free, we have $\beta(G) \le \frac{m}{2}-c_1m^{3/4 + q}$. Further, for every $m$ there is a digraph $G_m$ with $m$ edges such that $\overline{G_m}$ is $\mcB$-free and $\beta(G_m) \ge \frac{m}{2}-c_2m^{3/4 + q}$. 
\end{cor}

The above argument can be adapted to give a more general lower bound on $\beta(G)$ for digraphs $G$ that are not necessarily $\mcB$-free.
\begin{prop}
There is a constant $C > 0$ such that the following holds. For any undirected graph $H$ with $n$ vertices and $m$ edges, there is a digraph $G$ obtained by orienting the edges of $H$ such that 
$$\beta(G) \ge \frac{m}{2} -  C \log^{3/2}(n) \sqrt{mn}.$$ 
\end{prop}
\begin{proof}
We largely follow the argument and notation in the proof of Theorem~\ref{t:bfreelb}, omitting details that are identical. In particular, by adding isolated vertices, we may assume that $n$ is a power of two. Let $G$ be the orientation of $H$ that satisfies Lemma~\ref{lem:chernorient}. By our choice of orientation, we observe that
\begin{align*}
\pi(G)&\le \sum_{i = 0}^{\log_2(n/2)} \sum_{\substack{(A, B)  \in F_i}} 3\sqrt{\overline{e}(A, B)} \sqrt{2^i \log(en/2^i)}.
\end{align*}
Notice that the square root function is concave and thus for fixed $i$, the inner sum is maximized by setting all $\overline{e}(A, B)$ equal, i.e. setting all $\overline{e}(A, B) = \frac{m_i}{n \cdot 2^{-i-1}}$. This gives the desired upper bound on the directed surplus.
\begin{align*}
\pi(G) &\le 2 \sum_{i = 0}^{\log_2(n/2)} \sqrt{n \cdot 2^{-i-1}\cdot m_i\cdot 2^i \log(en/2^i)} \\
&\lesssim \sqrt{n} \sum_{i = 0}^{\log_2(n/2)} \sqrt{m_i \log(en/2^i)} \\
&\lesssim \sqrt{n \log n} \sum_{i = 0}^{\log_2(n/2)} \sqrt{m_i} \\
&\lesssim\log^{3/2}(n) \sqrt{mn} .
\end{align*}
\end{proof}

\section{Large girth and $r$-free digraphs}

When considering feedback arc sets in graphs with a forbidden subgraph, a natural choice for the forbidden subgraph is a cycle. Denote by $C_r^*$ the directed cycle on $r$ vertices and by $C_r$ the undirected cycle on $r$ vertices.
In an extreme case, if $G$ is free of all directed cycles of length up to $|V(G)|$ (i.e. $G$ is a \textit{directed acyclic graph}), then $\beta(G) = 0$.
We begin by discussing consider what the results from Section~\ref{s:boundbfree} give when the underlying undirected graph of $G$ is free from short cycles.
We first consider digraphs $G = (V, E)$ where we impose substantial structure on $G$ by requiring that the underlying undirected graph $\overline{G}$ has large \textit{girth}.

\begin{defn}
An undirected graph $H$ has \textit{girth} larger than $r$ if $H$ does not contain any cycle of length at most $r$, i.e. if $H$ is $\mcH$-free where $\mcH = \{C_3, \ldots C_{r}\}$.
\end{defn}

It follows from a result of Alon, Hoory, and Linial~\cite{AL02} that if an undirected graph $H$ has girth larger than $r$ for even $r \ge 4$, then $e(H) < v(H)^{1 + 2/r}$.
Combining this result with Theorem~\ref{t:bfreeub} gives an upper bound on $\beta(G)$ for digraphs whose underlying undirected graph has large girth.

\begin{cor}\label{p:undgirth} 
Let $r \ge 4$ be an even integer. There is a constant $c > 0$ such that if $G$ is a digraph with $m$ edges whose underlying undirected graph has girth larger than $r$, then, 
$$\beta(G) \le \frac{m}{2} - c m^{\frac{r+1}{r+2}}.$$
\end{cor}

Corollary \ref{p:undgirth} also holds if we only assume that the underlying undirected graph has no cycle of length $r$, except the constant $c$ depends on $r$. The proof uses the extremal bound for even cycles (see \cite{He}). This result is best possible up to the constant factor if for fixed even $r \geq 4$ we have $\textrm{ex}(n,C_{r})=\Theta(n^{1+2/r})$, as conjectured. 

We can hope to improve our bounds on $\beta(G)$ even when we only forbid specific oriented subgraphs.
\begin{defn}
A digraph $G = (V, E)$ is \textit{$r$-free} if it does not contain any $C_i^*$ with $i \leq r$ as a subgraph.
\end{defn}

Motivated by Corollary~\ref{p:undgirth} and the case of $3$-free digraphs, we pose the following stronger conjecture.
\begin{conj}
For any even integer $r \geq 4$ there is $c=c(r)>0$ such that if $G$ is an $r$-free digraph on $m$ edges and $n$ vertices, then  
$$\beta(G) \le \frac{m}{2} - c(r) m^{\frac{r+1}{r+2}}.$$
\end{conj}

If, as conjectured by Erd\H{o}s and Simonovits (c.f. Conjecture 4.10 in~\cite{FS13})
that $\ext(n, \mcB) = \Theta(n^{1 + 2/r})$ for $\mcB = \{C_3, \ldots, C_r\}$ for $r$ even, the above conjectured bound would be tight up to a constant function of $r$ in the surplus by Theorem~\ref{t:bfreelb}.

Although $\beta(G) = 0$ for any $n$-free digraph on $n$ vertices, $\beta(G)$ can remain quite large for $r$-free digraphs even when $r$ is relatively large as a function of $n$.
\begin{prop}\label{p:exrfree}
For any integer $r \geq 3$, if $n$ is a multiple of $r+1$, there is an $r$-free digraph $G$ on $n$ vertices such that $$\beta(G) = \frac{n^2}{(r + 1)^2}.$$
\end{prop}
\begin{proof}
Let $G$ be the balanced blowup of a directed cycle of length $r + 1$, obtained by taking the lexicographic product of $C_{r+1}^*$ with an empty graph on $t:=n/(r + 1)$ vertices. Digraph $G$ is $r$-free with shortest directed cycle of length $r+1$. 
The number of copies of $C_{r+1}^*$ is $t^{r+1}$, while each edge of $G$ is in exactly exactly $t^{r-1}$ copies of $C_{r+1}^*$. Thus, to remove all copies of $C_{r+1}^*$, we must remove at least $t^{r+1} /t^{r-1}= \frac{n^2}{(r + 1)^2}$ edges from $G$. If we remove the  $\frac{n^2}{(r + 1)^2}$ edges between consecutive parts in the blowup of $C_{r+1}^*$, the remaining directed subgraph is acyclic. Hence, $\beta(G) = \frac{n^2}{(r + 1)^2}$, as desired. 
\end{proof}

We suspect that the above construction is essentially optimal for $\beta(G)$ if $r$ is linear in $n$. We also believe the following more precise result when $r$ is very large. 

\begin{conj}\label{t:rfree}
If  $G = (V, E)$ is an $(r-1)$-free digraph with $r > 2n/3$, then $\beta(G) \le 1$.
\end{conj}

\begin{rem}
Conjecture~\ref{t:rfree} is tight as realized by the $(r-1)$-free digraph $G = (V, E)$ on $n = 3N$ vertices with vertex set $V = \{u_1, v_1, w_1, u_2, v_2, w_2 \ldots u_{N}, v_{N}, w_N\}$ (here indices are elements of $\mathbb{Z}_N$), with edges $(u_i, w_i), (v_i, w_i)$, $(w_i, v_{i+1})$, and $(w_i, u_{i+1})$ for $i \in \mathbb{Z}_N$. Then $r = 2N = \frac{2}{3}n$, but $\beta(G) = 2$.
\end{rem}

\section{Quasirandom directions}\label{s:qrandom} 
In this section, building on previous work of Griffiths~\cite{GRI13}, we give several equivalent characterizations of quasirandom directions of dense digraphs. One new characterization is that dense oriented graphs are quasirandom if and only if they have large minimum feedback arc set; in other words, a dense digraph $G$ on $n$ vertices with $m=\Omega(n^2)$ edges is quasirandom if and only if $\beta(G)=m/2 - o(m)$.  

It is possible to remove the assumption that $m= \Omega(n^2)$ by changing the property statements slightly, although the notions of quasirandomness we study would trivially hold when $m=o(n^2)$. Consequently, we can safely assume that the digraphs we consider are dense. We also use the asymptotic $o(\cdot)$ notation loosely. If we have two properties with $o(1)$ notation, so $P = P(o(1)), Q = Q(o(1))$, then $P$ implies $Q$ means that for each $\eps > 0$, there exists $\delta > 0$ so that if $G$ satisfies $Q(\delta)$ then it also satisfies $P(\eps)$.

Before stating our main result in this section, Theorem~\ref{t:expqrandom}, we first need some definitions. 


\begin{defn}
For a digraph $G = (V, E)$, call a $k$-tuple of vertices $(v_1, \ldots, v_k)$ an \textit{even-switch $k$-cycle} if the following conditions hold. For $i = 1, \ldots k$ (letting $v_{k+1} := v_1$) exactly one of $(v_i, v_{i+1}), (v_{i+1}, v_i) \in E$ and further, $(v_{i+1}, v_i) \in E$ for an even number of $i$. Let $E_k(G)$, be the number of distinct even-switch $k$-cycles in $G$ with respect to a labeling of $V$. Analogously, we define an \textit{odd-switch $k$-cycle}, and let $O_k(G)$ be the number of labeled odd-switch $k$-cycles in $G$.
\end{defn}

\begin{defn}
For a digraph $G = (V, E)$ on $n$ vertices, its adjacency matrix $A=A(G)$ is an $n \times n$ adjacency matrix with rows and columns indexed by vertices so that 
$$A_{uv} = \begin{cases}
1 & (u, v) \in E \\
-1 & (v, u) \in E \\
0 & \text{else.}
\end{cases}
$$
\end{defn}

\begin{defn}
For digraph $G = (V, E)$, we let $\bias_{\delta}(G)$ for $\delta \in (0, 1)$ be the size of the largest $\delta$-\textit{biased} subgraph of $G$, given by $$\bias_{\delta}(G) = \max \{ e(A, B) \mid A, B \subset V : e(B, A) \le \delta e(A, B)\}.$$
\end{defn}

We state an expanded version Theorem~\ref{t:qrandom} below.

\begin{thm}\label{t:expqrandom}
For a digraph $G$ on $n$ vertices and $m = \Omega(n^2)$ edges with underlying undirected graph $\overline{G}$, the following are equivalent.
\begin{enumerate}
\item $\tau(G) = o(m)$.
\item $\tau^*(G) = o(m)$.
\item $\pi(G) = o(m)$.
\item $N(C_4^{\rightarrow}, G) = \left( \frac18 + o(1) \right) N(C_4, \overline{G})$, where $C_4^{\rightarrow}$ is the orientation of a $C_4$ that has a homomorphism to an oriented edge.
\item For any labeling $L$ of $V$, $N_L(B, G) = \left( 2^{-|E(B)|} + o(1) \right) N_L(\overline{B}, \overline{G})$ for a fixed bipartite digraph $B$.
\item For any even $k \ge 4$, $E_k(G) = \left( \frac12+ o(1) \right) N_L(C_k, \overline{G})$.
\item For any even $k \ge 4$, $\Tr(A(G)^k) = o(\Tr(A(\overline{G})^k)))$.
\item $|\lambda_1(G)| = o(|\lambda_1(\overline{G})|)$.
\item  There exists some $\delta = 1 - o(1)$ such that $\bias_{\delta}(G) \le o(m)$.
\end{enumerate}
If any of these equivalent conditions is satisfied, we say that $G$ has \emph{quasirandom direction} with respect to the underlying undirected graph $\overline{G}$.
\end{thm}
Many of the above implications are found in~\cite{GRI13} or arise from variants of earlier work such as~\cite{FAN91}. Our primary contribution is to relate the size of a digraph's minimum feedback arc set (via $\pi(G)$) to characterizing quasirandom direction. We highlight a few additional interesting implications as well. 

It is natural to consider the property that a digraph is \textit{almost balanced}, i.e.~almost all vertices have nearly equal indegree and outdegree. This turns out to be insufficient for guaranteeing quasirandomness (as seen by the balanced blowup of a directed triangle). However, when the indegree and outdegree of graph vertices are similar, we do find that edges in a graph are balanced with respect to partitions of the graph. More precisely, we have the following result.

\begin{prop}\label{prop:identity}
For every digraph $G=(V,E)$, we have $$\sum_{v \in V} |d^+(v) - d^-(v)|=2\tau_{\sqcup}(G).$$
\end{prop}
\begin{proof}
Consider a partition $V=A \sqcup B$. We have 
\begin{eqnarray*}\sum_{v \in A} d^+(v) - d^-(v) &= & e(A,V)-e(V,A) \\ & = & e(A,B)+e(A,A)-(e(B,A)+e(A,A)) \\ & = & e(A,B)-e(B,A).\end{eqnarray*}
Note that the above expression is maximized by taking $A$ to be the set of all vertices $v \in V$ such that $d^+(v) > d^-(v)$, and the expression in this case is equal to $\tau_{\sqcup}(G)$. For this particular set $A$, we have $$\sum_{v \in A} d^+(v) - d^-(v)=\frac{1}{2}\sum_{v \in V} |d^+(v) - d^-(v)|,$$ as the sum of the indegrees equals the sum of the outdegrees of vertices in a graph. Putting this together gives the desired identity.  
\end{proof}
Since $\tau_{\sqcup}(G) \le \tau(G)$, Proposition~\ref{prop:identity} implies that being almost balanced is a \textit{necessary} condition for a digraph to be quasirandomly oriented.

We will show the implications of Theorem~\ref{t:qrandom} in a series of lemmas that follow, leveraging the results of~\cite{GRI13,KAL13} that show some of the equivalence directions.

\begin{lemma}\label{l:pitau}
Every digraph $G$ with $n$ vertices satisfies $$\pi(G) \le n\sqrt{\tau^*(G)}.$$
\end{lemma}

\begin{proof}
Fix an optimal ordering $\rho$ of $V(G)$ with $\beta(G)$ backwards edges. With foresight, let $t=\lfloor n\tau^*(G)^{-1/2}\rfloor$. Note that $t \geq 2$ as $\tau^*(G) \leq n^2/4$ follows from the definition of $\tau^*$. Fix an equitable partition of the vertices in $V$ into $t$ subsets $A_{1},...,A_{t}$,  such that for $i<j$, the vertices in $A_{i}$ come strictly before the vertices in $A_{j}$ with respect to $\rho$. 
The number of unordered pairs of vertices with both vertices in the same part is $\sum_{i=1}^t {|A_i| \choose 2} \leq n^2/(2t)$, where the upper bound follows from the partition being equitable. Let $B_{i}=\bigcup_{j=i+1}^{t}A_{j}$. Noting that $A_i$ and $B_i$ are disjoint, we see that 
\begin{equation}\label{tau*sub}\tau^{*}(G)\geq e(A_{i},B_{i})-e(B_{i},A_{i}).\end{equation}
Observe that $ \sum_{i=1}^{t}e(A_{i},B_{i})$ counts the number of forward edges in $G$ with vertices in different parts by ordering $\rho$. 
By our optimal choice of $\rho$, there are $\frac{m}{2}+\pi(G)$ forwards edges, and therefore 
\begin{equation}\label{backineq} \sum_{i=1}^{t}e(A_{i},B_{i}) +\frac{n^2}{2t} \geq\frac{m}{2}+\pi(G).\end{equation}
So $\beta(G)$, which is the number of backwards edges with respect to $\rho$, is at least 
\begin{align*}
\sum_{i=1}^{t}e(B_{i},A_{i})
&\ge \frac{m}{2} + \pi(G) -  \sum_{i=1}^{t}(e(A_{i},B_{i})-e(B_{i},A_{i})) - \frac{n^2}{2t} \\
&\geq \frac{m}{2} + \pi(G) - t\tau^{*}(G)-\frac{n^2}{2t} \\
&\geq \frac{m}{2} + \pi(G) -2n\sqrt{\tau^{*}(G)},
\end{align*}
where the first inequality is by Inequality (\ref{backineq}), the second inequality substitutes Inequality (\ref{tau*sub}) for each $i$, and the last 
inequality utilizes the choice of $t$. This implies that $\pi(G) = m/2 - \beta(G) \le n \sqrt{\tau^*(G)}$.
\end{proof}

\begin{lemma}\label{l:ektr}
For a dense digraph $G = (V, E)$ with underlying undirected graph $\overline{G}$ and some labeling $L$ of $V$, for any even integer $k \ge 4$, 
$$\Tr(A^k) = 2 E_k(G) - \Tr(A(\overline{G})^k),$$
where $A(\overline{G})$ is the adjacency graph of $\overline{G}$, defined as usual.
Thus, $$E_k(G) = \left(\frac12 + o(1)\right) N_L(C_k, \overline{G}) \text{ if and only if } \Tr(A^k) = o\left(\Tr(A(\overline{G})^k)\right).$$
\end{lemma}
\begin{proof}
The proof extends the argument of Kalyanasundaram and Shapira (Claim 2.2 in~\cite{KAL13}) who proved an analogous result for tournaments. The $(v, v)$ entry of $A^k$ is the number of even-switch $k$-cycles with vertex $v$ minus the number of odd-switch $k$-cycles with vertex $v$. Thus, $\Tr(A^k) = E_k(G) - O_k(G)$. Note that $E_k(G) + O_k(G)  = \Tr\left(A(\overline{G})^k\right).$
This gives  
that 
$$E_k(G) = \frac12 \left( \Tr(A^k) +  \Tr(A(\overline{G})^k)\right).$$
Observe that $\Tr\left(A(\overline{G})^k\right)$ counts the number of closed walks of length $k$ in $\overline{G}$, which for dense $\overline{G}$ is asymptotically equal to the number of labelled copies of $C_k$ in $\overline{G}$. 
Consequently, $E_k(G) = (\frac12 + o(1)) N_L(C_k, \overline{G})$ if and only if $\Tr(A^k) = o(\Tr(A(\overline{G})^k))$. 
\end{proof}

\begin{proof}[Proof of Theorem~\ref{t:expqrandom}]
Lemma~\ref{l:tautau} implies that (1) is equivalent to (2). The proof of 
Lemma~\ref{t:betatau} shows that $\frac12\tau^*(G) \le \pi(G)$, so $(3)$ implies $(1)$. Lemma~\ref{l:pitau} shows that (2) implies (3) as $m = \Omega(n^2)$, so if $\tau^*(G) = o(m)$, then $\pi(G) = O(n \sqrt{\tau^*(G)}) = n \cdot o(m^{1/2}) = o(m)$. 
By setting $B = C_4$ and considering all possible valid labelings of a directed $C_4$ as a subgraph of $G$, we see that (5) implies (4) and similarly, (5) implies (6). Theorem 1.1 of~\cite{GRI13} shows that (1) is equivalent to (5) and that (4) and (6) each imply (1), which gives that (1),(4),(5), (6) are all equivalent. Lemma~\ref{l:ektr} shows that (6) is equivalent to (7). As noted in~\cite{KAL13}, (7) is equivalent to (8) (and they also show that (7) is equivalent to (1) directly). We finally observe that (9) is equivalent to (1) by our assumption that $m = \Theta(n^2)$, so $\tau(G) = o(n^2)$ exactly if there exists some $\delta = 1 - o(1)$ such that $\bias_{\delta }(G) = o(n^2) = o(m)$.
\end{proof}

\appendix

\section{Simplified proof of Theorem~\ref{t:betaub}}\label{a:betaub}

Below we give a simplified, but similarly motivated proof of Theorem~\ref{t:betaub}, originally shown in~\cite{BS97}, also observing the bound $\beta(G) = \frac{m}{2} - \Omega(m^{3/4})$ that was not stated in the earlier work but that follows immediately from the bound $\beta(G) \le \frac{m}{2} - c \sum_{v \in V} \sqrt{d(v)}$.

\begin{rem}
In order to construct a small feedback arc set of digraph $G$ with $n$ vertices, a natural idea is to take a greedy approach. More precisely, given an ordering $v_1,\ldots,v_n$ of the vertices of $G$, consider the following greedy algorithm for constructing another ordering $u_1,\ldots,u_n$ for which the set of directed edges $(u_j,u_i)$ with $j>i$ is small (this yields a small feedback arc set). At step $i$, we will fix the relative order of $v_1,\ldots,v_i$, for now as $u_1, \ldots u_i$. We place $v_{i+1}$ in the sequence so as to minimize the number of back edges, meaning that over all 
$t \in \{0,1,\ldots,i\}$,  we minimize the number of edges either of the form $(v_{i+1},u_j)$ with $j \leq t$ or of the form $(u_j,v_{i+1})$ with $t < j \leq i$. We update the ordering by keeping $u_j$ the same for $j \leq t$, letting $u_{t+1}=v_{i+1}$, and letting $u_{j+1}$ be the previous $u_j$ for $t < j \leq i$.

While the exhaustive search for finding a minimum feedback arc set over all $n!$ orderings of an $n$-vertex graph is slow, and determining the size of a minimum feedback arc set is known to be $\NP$-hard, given the ordering $v_1,\ldots,v_n$ of the vertices of a digraph $G$, the greedy algorithm runs in polynomial-time. 

The above greedy algorithm is somewhat challenging to analyze and it is unclear how to find a good ordering of the vertices of $G$ to start with. Consequently, in the following proof of Theorem~\ref{t:betaub}, we consider a restricted version of the above procedure applied to a random ordering of the vertex set. Such an algorithm is much simpler to analyze, and we are thus able to use it to provably construct a small feedback arc set of any digraph $G$.
\end{rem}

\begin{proof}[Proof of Theorem~\ref{t:betaub}]
Let $G$ be  a digraph with $n$ vertices. 
A restricted greedy algorithm for finding a feedback arc set is as follows. Again, we are given a digraph $G$ with $n$ vertices and an ordering $v_1,\ldots,v_n$ of the vertices of $G$. At step $i$, we again fix the relative order of $v_1,\ldots,v_i$. However, instead of placing $v_{i+1}$ optimally given the relative order of $v_1,\ldots,v_i$, we place $v_{i+1}$ before or after all of the vertices $v_1,\ldots,v_i$, whichever minimizes the number of backwards edges. In this restricted greedy algorithm, letting $S_i=\{v_1,\ldots,v_i\}$, we obtain a feedback arc set of size $\sum_{i=0}^{n-1} \min\left(d^+(v_{i+1},S_i),d^-(v_{i+1},S_i)\right)$. In particular, from the restricted greedy algorithm, we get the following bound on the directed surplus:  \begin{equation}\label{eq123}\pi(G) \geq \frac{1}{2}\sum_{i=0}^{n-1} \left|d^+(v_{i+1},S_i)-d^-(v_{i+1},S_i)\right|\end{equation} for every ordering $v_1,\ldots,v_n$ of the vertices of $G$. We observe that given an ordering, the restricted greedy algorithm takes time $O(n^2)$ as each of the $n$ steps takes time $O(n)$.

We pick the starting ordering uniformly at random. A given vertex $v$ will appear as a $v_{i+1}$ with $i \in \{0,1,\ldots,n-1\}$ uniformly. 
The expected value of the contribution $|d^+(v,S_i)-d^-(v,S_i)|$ vertex $v$ gives when chosen $v_{i+1}$ is equal to the expected value of the following random variable $X_v$, defined as follows.
For a vertex $v$, randomly order $\{v\} \cup N^+(v) \cup N^-(v)$, and let $X_v$ be the absolute difference of the number of in-neighbors and out-neighbors that come before $v$ in this ordering. In other words, the distribution of $X_v$ is that of the absolute value of the sum of $h$ random numbers taken without replacement from a multiset of $1$'s with multiplicity $d^+(v)$ and $-1$'s with multiplicity $d^{-}(v)$, where $h$ is a random integer picked uniformly from $0$ to $d(v)$. Similar to the fact that a binomial random variable $B(d,1/2)$ or hypergeometric random variable has standard deviation on the order of $\sqrt{d}$, we get that the expected contribution from $v$ is $\Omega\left(\sqrt{d(v)}\right)$. In total, we get $\pi(G)=\Omega\left(\sum_{v \in V(G)}\sqrt{d(v)}\right)$ as desired. The last desired inequality  follows from the next lemma. 
\end{proof}

\begin{lem}
If $G=(V,E)$ is a graph with $m$ edges, then $\sum_{v \in V}\sqrt{d(v)} \geq \frac{1}{4}m^{3/4}$. 
\end{lem}
\begin{proof}
At most $m/2$ edges of $G$ are contained within the subgraph of $G$ induced by the $\sqrt{m}$ vertices of highest degree. Further, there are fewer than $\sqrt{m}$ vertices of degree larger than $2 \sqrt{m}$. Thus, 
there are more than $m/2$ edges that are incident to a vertex of degree at most $2\sqrt{m}$, and hence 
$$\sum_{v \in V} \sqrt{d(v)} \ge \sum_{v : d(v) \le 2\sqrt{m}} \sqrt{d(v)} = \sum_{v : d(v) \le 2\sqrt{m}}  \frac{d(v)}{\sqrt{d(v)}} \ge \sum_{v : d(v) \le 2\sqrt{m}}  \frac{d(v)}{2 m^{1/4}} \ge \frac14 m^{3/4}. $$
\end{proof}

\end{document}